\documentclass{article}
\usepackage{amssymb,amsmath,amsthm}
\newtheorem{theorem}{Theorem}

\newtheorem{lemma}{Lemma}

\newtheorem{remark}{Remark}
\date{}
\numberwithin{equation}{section} \numberwithin{theorem}{section}
\numberwithin{lemma}{section} \numberwithin{corollary}{section}
\numberwithin{remark}{section} \numberwithin{proposition}{section}
\numberwithin{definition}{section}

\def\Xint#1{\mathchoice
{\XXint\displaystyle\textstyle{#1}}%
{\XXint\textstyle\scriptstyle{#1}}%
{\XXint\scriptstyle\scriptscriptstyle{#1}}%
{\XXint\scriptscriptstyle\scriptscriptstyle{#1}}%
\!\int}
\def\XXint#1#2#3{{\setbox0=\hbox{$#1{#2#3}{\int}$ }
\vcenter{\hbox{$#2#3$ }}\kern-.6\wd0}}

\def\dashint{\Xint-}

\begin{document}
\newcommand{\n}{\noindent}
\newcommand{\vs}{\vskip}
\newcommand{\avint}{\mathop{\ooalign{$\int$\cr$-$}}}

\title{H\"{o}lder continuity for the solutions of the p(x)-Laplace equation
with general right-hand side}

\author{A. Lyaghfouri \\
United Arab Emirates University\\
Department of Mathematical Sciences\\
Al Ain, Abu Dhabi, UAE} \maketitle

\begin{abstract}
We show that bounded solutions of the quasilinear elliptic equation
$\Delta_{p(x)} u=g+div(\textbf{F})$ are locally H\"{o}lder continuous
provided that the functions $g$ and $\textbf{F}$ are in suitable
Lebesgue spaces.

\end{abstract}


\n Key words : $p(x)$-Laplacian, H\"{o}lder continuity.

\vs 0.3cm

\n MSC: 35B65, 35J92

\section{Introduction}\label{S:intro}

\n We consider the following equation

\begin{equation}\label{1.1}
\Delta_{p(x)} u=g+div(\textbf{F})\quad\text{in } W^{-1,q(x)}(\Omega)
\end{equation}

\n where $\Delta_{p(x)} u=div(|\nabla u|^{p(x)-2}\nabla u)$ is the
$p(x)$-Laplacian, $\Omega$ is an open bounded domain of
$\mathbb{R}^n$, $n\geq 2$, $x=(x_1,...,x_n)$,
$q(x)={{p(x)}\over{p(x)-1}}$, and $p:~\Omega \rightarrow \mathbb{R}$ is a
measurable function which satisfies for some positive constants $p_{+}>p_{-}>1$
and $s>n$
\begin{eqnarray}\label{1.2-1.3}
&&p_{-}  \leq p(x) \leq p_{+}\quad\text{a.e. } x\in \Omega\\
&&\nabla p\in \big(L^s(\Omega)\big)^n.
\end{eqnarray}

\n As a consequence of (1.2) and (1.3), we have $p\in W^{1,s}(\Omega)$. Moreover, due to
Sobolev embedding $W^{1,s}(\Omega)\subset C^{0,\beta}(\Omega)$,
$\displaystyle{\Big(\beta=1-{n\over{s}}\Big) }$,
$p$ is H\"{o}lder continuous in $\Omega$.

\n We call a solution of equation (1.1) any function $u\in W^{1,p(x)}(\Omega)$
that fulfills
\[\int_{\Omega} |\nabla u|^{p(x)-2}\nabla u.\nabla\zeta dx=-\int_{\Omega} g(x)\zeta dx+\int_{\Omega} \textbf{F}(x).\nabla\zeta dx\quad
\forall\zeta\in W_0^{1,p(x)}(\Omega)\]
\n The Lebesgue and Sobolev spaces with variable exponents are defined (see for example \cite{[AS]},
\cite{[FZ]} and \cite{[KR]}) by:

\begin{eqnarray*}
&&L^{p(x)} (\Omega) = \Big\{u : \Omega \rightarrow \mathbb{R} ~\mbox{ measurable}~ :~\int_{\Omega}
|u(x)|^{p(x)}dx <  \infty~\Big\}\\
&&W^{1,p(x)} (\Omega) = \Big\{~u \in L^{p(x)}\Omega)~ /~ \nabla u \in \big(L^{p(x)}(\Omega)\big)^n~\Big\}\\
&&W_0^{1,p(x)}(\Omega)=\overline{C_0^{\infty}(\Omega)}_{W^{1,p(x)}(\Omega)}
\end{eqnarray*}

\n These spaces are separable, complete and reflexive, when equipped with the following norms
\begin{eqnarray*}
&&\|u\|_{L^{p(x)} (\Omega)} = \inf\left\{\lambda > 0 ~/~ \int_{\Omega}
\left|\frac{u}{\lambda}\right|^{p(x)}dx \leq 1~\right\}\\
&&\|u\|_{1,p(x)} = \|u\|_{p(x)} + \|\nabla u\|_{p(x)},\quad \|\nabla u\|_{p(x)} = \displaystyle{\sum^{n}_{i=1}
\left\|\frac{\partial u}{\partial x_{i}}\right\|_{p(x)}}
\end{eqnarray*}

\n Our aim is to establish H\"{o}lder continuity for bounded solutions
of (1.1). We observe that if $p(x)>n$ in an open set $U\subset\subset\Omega$,
then by Sobolev embedding $W^{1,p(x)}(U)\subset W^{1,p_m(U)}(U)\subset C^{0,\alpha}(U)$,
where $\displaystyle{p_m(U)=\min_{x\in \overline{U}} p(x)}$ and $\displaystyle{\alpha=1-{n\over{p_m(U)}} }$.
Therefore any solution of (1.1) is H\"{o}lder continuous in $U$.
In this paper, we assume that
\begin{eqnarray}\label{1.4}
p(x)\leq n\quad \forall x\in \Omega.
\end{eqnarray}

\vs 0,3cm\n $g$ is a real valued function that satisfies
for a positive number $t_1$
\begin{eqnarray}\label{1.5-1.6}
&&t_1>{{n}\over{p(x)}}\quad \forall x\in \Omega\\
&& g\in L^{t_1}(\Omega).
\end{eqnarray}

\vs 0,3cm\n $\textbf{F}=(F_1,...,F_n)$ is a vector function that satisfies
for a positive number $t_2$
\begin{eqnarray}\label{1.7-1.8}
&&t_2>{{n}\over{p(x)-1}}\quad \forall x\in \Omega\\
&& \textbf{F}\in L^{t_2}(\Omega).
\end{eqnarray}

\n Among problems that fit in the equation (1.1) setting,
is the dam problem ($g=0$ and $p(x)$ constant, \cite{[CaL]}, \cite{[CL1]}, \cite{[CL2]}).
Another problem is the obstacle problem ($\textbf{F}=0$, \cite{[CL6]}, \cite{[CL7]},
\cite{[CLR]}, \cite{[CLRT]}, \cite{[RT]}).
For a more general framework, we refer to \cite{[CL4]}, \cite{[CL8]} and \cite{[L5]}.

\vs 0.2cm\n The main result of this paper is the following theorem.
\begin{theorem}\label{t1.1} Assume that (1.2)-(1.8) hold. Then any bounded solution
of (1.1) is locally H\"{o}der continuous in $\Omega$.
\end{theorem}

\vs 0.2cm\n When the righthand side of (1.1) is a nonnegative Radon measure with a
suitable growth condition, the regularity result in Theorem 1.1 was established in
\cite{[K]} for the $p-$laplacian, and extended in \cite{[CL3]} and \cite{[L4]}
respectively for the $A-$laplacian and $p(x)-$laplacian.

\vs 0.2cm\n Throughout this paper, we will denote by $B_r(x)$ (resp. $\overline{B}_r(x)$) 
the open (resp. closed) ball of center $x$ and radius $r$  in $\mathbb{R}^n$, with
$\omega_n=|B_1|$ standing for the measure of the unit open ball $B_1$.

\section{Proof of Theorem 1.1}\label{2}

\vs 0.3cm\n In this section, we denote by $u$ a bounded solution of (1.1)
and will show that it is locally H\"{o}lder continuous in $\Omega$.
The proof is based on Lemma 2.1.

\begin{lemma}\label{l2.1}
Let $M=||u||_\infty$, $\delta={{s-n}\over{2sn}}$, $q={{n+s}\over 2}$, $m={{sq(1+\delta)}\over{s-q}}$, and
\[\displaystyle{R_1=\min\Big(\text{diam}(\Omega)/16,c(n,s,p_-,p_+)
\Big( \int_{\Omega} |\nabla u|^{p(x)}dx+1\Big)^{{-m}\over{\beta}}}\Big).\]
Assume that (1.2)-(1.8) hold. Then there exists a positive constant
\[C=C(n,s,p_-,p_+,M,\text{diam}(\Omega),||\nabla u||_{L^{p(x)}(\Omega)}, |\nabla p|_{L^s(\Omega)}, ||g||_{L^{t_1}(\Omega)}, ||F||_{L^{t_2}(\Omega)} ),\]
such that we have for $R\in(0, R_1)$ with $B_{2R}(x_0)\subset \Omega$
\begin{equation}\label{e2.1}
\int_{B_r(x_0)} |\nabla u|^{p_m}dx \leq  C r^{n-p_m+\alpha
p_m}\qquad \forall r\in (0, R),
\end{equation}
where $\displaystyle{p_m=\min_{x\in \overline{B}_R(x_0)}p(x)}$ and
$\displaystyle{\alpha=1-\max\left({{n}\over{t_1p_m}},{{n}\over{t_2(p_m-1)}}\right) }$.
\end{lemma}

\vs 0,3cm\n The proof of Lemma 2.1 is based on Lemma 2.2.

\begin{lemma}\label{l2.2}
Assume that (1.2)-(1.8) hold and let $R_1$ be the positive number in Lemma 2.1.
Then there exists a positive constant
\[C_1=C_1(n,s,p_-,p_+,\text{diam}(\Omega),||\nabla u||_{L^{p(x)}(\Omega)}, ||\nabla p||_{L^s(\Omega)} ),\]
such that we have for $R\in(0,R_1)$, $B _{2R}(x_0)\subset \Omega$, $v\in W^{1,p(x)} (B_R(x_0))$
with $\Delta_{p(x)}v = 0$ in $B_R(x_0)$ and $v=u$ on $\partial B_{R}(x_0)$
\begin{equation}\label{e2.2}
\int_{B_r(x_0)} |\nabla v|^{p(x)}dx \leq C_1\Big(\Big( \Big({ r\over R}\Big)^n +R^{{\beta}\over{2}}\Big)\int_{B_R(x_0)} |\nabla
u|^{p(x)}dx + R^{n+\beta}\Big)\qquad \forall r\in (0, R).
\end{equation}
\end{lemma}

\n {\bf Proof.} First, observe that we have
$v\in C_{loc}^{1,\sigma}(B_R(x_0))$ (see for example \cite{[CM]}). Next, since
$\displaystyle{\int_{B_R(x_0)} {{|\nabla v|^{p(x)}}\over{p(x)}} dx}$ minimizes the energy
$\displaystyle{\int_{B_{R}(x_0)} {{|\nabla \xi|^{p(x)}}\over {p(x)}} dx}$ over all
functions $\xi\in W^{1,p(x)}( B_{R}(x_0))$ such that $\xi=u$ on $\partial B_{R}(x_0)$, we get
\begin{eqnarray}\label{e2.3}
\int_{B_R(x_0)} |\nabla v|^{p(x)}dx &=& \int_{B_R(x_0)}p(x) {{|\nabla v|^{p(x)}}\over{p(x)}} dx~\leq~ p_+\int_{B_R(x_0)}{{|\nabla v|^{p(x)}}\over{p(x)}}dx\nonumber\\
&\leq& p_+\int_{B_R(x_0)}{{|\nabla u|^{p(x)}}\over{p(x)}}dx~\leq~{p_+\over p_-}\int_{B_R(x_0)} |\nabla u|^{p(x)}dx.
\end{eqnarray}

\n Using (2.3), we get for $r\in \Big[{R\over 8}, R\Big)$
\begin{eqnarray*}
\int_{B_r(x_0)} |\nabla v|^{p(x)}dx &=& \Big( { R\over
r}\Big)^n\Big( { r\over R}\Big)^n\int_{B_r(x_0)} |\nabla v|^{p(x)}dx
~\leq~ 8^n \Big( { r\over R}\Big)^n\int_{B_r(x_0)} |\nabla v|^{p(x)}dx \\
&\leq& 8^n \Big( { r\over R}\Big)^n\int_{B_R(x_0)} |\nabla v|^{p(x)}dx
~\leq~ 8^n {p_+\over p_-}\Big( { r\over R}\Big)^n\int_{B_R(x_0)} |\nabla
u|^{p(x)}dx.
\end{eqnarray*}

\n Therefore it is enough to prove (2.2) for $r\in \Big(0, {R\over 8}\Big)$.
Since $\Delta_{p(x)}v = 0$ in $B_R(x_0)$, $v=u$ on $\partial B_R(x_0)$, we
obtain by the maximum principle $||v||_{L^\infty ( B_R(x_0)) }\leq M$. Moreover, we know
(see \cite{[CL5]}, Corollary 2.1) for $\gamma=\beta/2$ and $\epsilon_0={1\over 2}\min(1,m-1)$,
that there exist two positive constants $\displaystyle{R_0=R_0\Big(n,s,p_-,p_+,\text{diam}(\Omega),\int_{B_R(x_0)} |\nabla v|^{p(x)}dx\Big)}$
and $c_1=c_1(n,s,p_-,p_+)$ such that for each $R\in(0,\min(R_0,\text{diam}(\Omega)/16))$, $r\in(0,R/8)$,
and $K=||\nabla p||_{L^s(\Omega)}\big(1+(\text{diam}(\Omega)/2)^\beta|\nabla p|_{L^s(\Omega)}\big)$, we have
\begin{eqnarray*}
\dashint_{B_r(x_0)} |\nabla
v|^{p(x)}dx &\leq&c_1\Big(\dashint_{B_R(x_0)} |\nabla v|^{p(x)}dx \nonumber\\
&&\quad+KR^\beta\Big(1+\int_{B_R(x_0)} |\nabla v|^{p(x)}dx\Big)^{\epsilon_0(1+\delta)}
\Big(1+\dashint_{B_R(x_0)} |\nabla v|^{p(x)}dx\Big)^{1+\delta}\Big)\nonumber
\end{eqnarray*}
\n where $\displaystyle{\dashint_E f dx={1\over{|E|}}\int_E f dx}$ denotes the mean 
value of the function $f$ on the measurable set $E$.

\vs0.1cm\n Using (2.3), we obtain for
$c_2=\displaystyle{c_1\Big({p_+\over p_-}\Big)^{(1+\epsilon_0)(1+\delta)}\max\Big(1,2^\delta K\Big(1+{p_+\over p_-}\int_{\Omega} |\nabla u|^{p(x)}dx\Big)^{\epsilon_0(1+\delta)}\Big)}$

\begin{eqnarray*}
\dashint_{B_r(x_0)} |\nabla
v|^{p(x)}dx &\leq&c_1\Big({p_+\over p_-}\dashint_{B_R(x_0)} |\nabla u|^{p(x)}dx \nonumber\\
&&\quad+KR^\beta\Big(1+{p_+\over p_-}\int_{B_R(x_0)} |\nabla u|^{p(x)}dx\Big)^{\epsilon_0(1+\delta)}
\Big(1+{p_+\over p_-}\dashint_{B_R(x_0)} |\nabla u|^{p(x)}dx\Big)^{1+\delta}\Big)\nonumber\\
&\leq& c_1\Big({p_+\over p_-}\dashint_{B_R(x_0)} |\nabla u|^{p(x)}dx \nonumber\\
&&\quad+KR^\beta\Big(1+{p_+\over p_-}\int_{\Omega} |\nabla u|^{p(x)}dx\Big)^{\epsilon_0(1+\delta)}
\Big(1+{p_+\over p_-}\dashint_{B_R(x_0)} |\nabla u|^{p(x)}dx\Big)^{1+\delta}\Big)\nonumber\\
&\leq& c_2\Big(\dashint_{B_R(x_0)} |\nabla u|^{p(x)}dx+R^\beta+R^\beta\Big(\dashint_{B_R(x_0)} |\nabla u|^{p(x)}dx\Big)^{1+\delta}\Big)\nonumber
\end{eqnarray*}

\n which leads to
\begin{eqnarray*}
\int_{B_r(x_0)} |\nabla v|^{p(x)}dx &\leq& c_2\Big(\Big({ r\over R}\Big)^n\int_{B_R(x_0)} |\nabla u|^{p(x)}dx+\omega_nr^nR^\beta\nonumber\\
&&\quad+\omega_n^{-\delta}r^nR^{\beta-n(1+\delta)}\Big(\int_{B_R(x_0)} |\nabla u|^{p(x)}dx\Big)^{1+\delta}\Big)\nonumber
\end{eqnarray*}

\n or
\begin{eqnarray}\label{e2.4}
\int_{B_r(x_0)} |\nabla v|^{p(x)}dx &\leq& c_2\Big(\Big({ r\over R}\Big)^n\int_{B_R(x_0)} |\nabla u|^{p(x)}dx+\omega_nR^{n+\beta}\nonumber\\
&&\quad+\omega_n^{-\delta}R^{\beta-n\delta}\Big(\int_{B_R(x_0)} |\nabla u|^{p(x)}dx\Big)^{1+\delta}\Big)\nonumber\\
&\leq&c_2\Big(\Big({ r\over R}\Big)^n\int_{B_R(x_0)} |\nabla u|^{p(x)}dx+\omega_nR^{n+\beta}\nonumber\\
&&\quad+\omega_n^{-\delta}R^{\beta\over 2}\Big(\int_{\Omega} |\nabla u|^{p(x)}dx\Big)^{\delta}\Big(\int_{B_R(x_0)} |\nabla u|^{p(x)}dx\Big)\Big).
\end{eqnarray}

\vs 0.1cm\n We point out (see \cite{[AM]}, \cite{[CL5]}) that
$$R_0\approx \Big( { \epsilon_0\over {c(n,p_-,p_+)^m}c(\beta)}\Big)^{2/\beta}
\Big( \int_{B_R(x_0)} |\nabla v|^{p(x)}dx+1\Big)^{{-2m\epsilon_0}\over{\beta}}$$
\n Using (2.3) again, we get for some positive constant $c_0=c_0(n,s,p_-,p_+)$ independent of $R$
\begin{eqnarray*}
R_0&\geq&c_0\Big( \int_{B_R(x_0)} |\nabla v|^{p(x)}dx+1\Big)^{{-2m\epsilon_0}\over{\beta}}
~\geq~c_0\Big( {{p_+}\over{p_-}}\int_{B_R(x_0)} |\nabla u|^{p(x)}dx+1\Big)^{{-2m\epsilon_0}\over{\beta}}\\
&\geq&c_0\Big({p_+\over p_-}\Big)^{{-2m\epsilon_0}\over{\beta}}\Big(\int_{\Omega} |\nabla
u|^{p(x)}dx+1\Big)^{{-m}\over{\beta}}
~=~c(n,s,p_-,p_+)\Big( \int_{\Omega} |\nabla u|^{p(x)}dx+1\Big)^{{-m}\over{\beta}}
\end{eqnarray*}

\n Thus, we can take $\displaystyle{R_1=\min\Big(\text{diam}(\Omega)/16,c(n,s,p_-,p_+)
\Big( \int_{\Omega} |\nabla u|^{p(x)}dx+1\Big)^{{-m}\over{\beta}}}\Big)$, and (2.2) follows
from (2.4) if we choose
\[\displaystyle{C_1=c_2\max\Big(1,\omega_n,\omega_n^{-\delta}\Big(\int_\Omega |\nabla u|^{p(x)}dx\Big)^{\delta}\Big)=C_1(n,s,p_-,p_+,\text{diam}(\Omega),||\nabla p||_{L^s(\Omega)},||\nabla u||_{L^{p(x)}(\Omega)})}.\]
\qed

\vs 0,3cm\n\emph{Proof of Lemma 2.1.} Let $R_1$ be as in Lemma 2.1, $v$ as in Lemma 2.2, and $R\in(0,R_1)$.
First, we have for $ r\in (0, R)$
\begin{align}\label{2.5}
   & \int_{B_r(x_0)} |\nabla u|^{p(x)} =\int_{B_r(x_0)} |\nabla u|^{p(x)-2} \nabla u .\nabla (u-v)dx+
   \int_{B_r(x_0)} |\nabla u|^{p(x)-2} \nabla u .\nabla v dx \nonumber\\
& = \int_{B_r(x_0)} \Big(|\nabla u|^{p(x)-2} \nabla u -
   |\nabla v|^{p(x)-2} \nabla v \Big).\nabla ( u-v)dx +
     \int_{B_r(x_0)} |\nabla v|^{p(x)-2} \nabla v .\nabla (u-v)dx \nonumber\\
     &+ \int_{B_r(x_0)} |\nabla u|^{p(x)-2} \nabla u .\nabla v dx
     = I_1 + I_2 + I_3.
\end{align}
\n Next, using the monotonicity of $\xi~\rightarrow~|\xi|^{p(x)-2}\xi$,
the fact that $\Delta_{p(x)}v = 0$ in $B_R(x_0)$ and $u=v$ on $\partial B_{R}(x_0)$, we obtain
\begin{eqnarray}\label{2.6}
  I_1&\leq &\int_{B_{R}(x_0)} \Big(|\nabla u|^{p(x)-2} \nabla u -
   |\nabla v|^{p(x)-2} \nabla v \Big).\nabla ( u-v)dx\nonumber\\
   &=& \int_{B_{R}(x_0)}|\nabla u|^{p(x)-2} \nabla u .\nabla ( u-v)dx - \int_{B_{R}(x_0)}
   |\nabla v|^{p(x)-2} \nabla v .\nabla ( u-v)dx\nonumber\\
   &=& -\int_{B_{R}(x_0)} g.(u-v)dx+\int_{B_{R}(x_0)} F.\nabla ( u-v)dx
\end{eqnarray}
\n Applying Young's inequality, we derive for some positive constant $c$ depending
only on $p_-$ and $p_+$
\begin{eqnarray}\label{2.7}
I_2&=&\int_{B_r(x_0)} |\nabla v|^{p(x)-2}\nabla v.\nabla u dx-\int_{B_r(x_0)} |\nabla
v|^{p(x)}~\leq ~ \int_{B_r(x_0)}   |\nabla v|^{p(x)-1}  .|\nabla u| dx\nonumber\\
&\leq&  {1\over 4}  \int_{B_r(x_0)}  |\nabla u|^{p(x)} dx
+ c\int_{B_r(x_0)}  |\nabla v|^{p(x)}dx.
\end{eqnarray}
\n Similarly, we obtain
\begin{align}\label{2.8}
   &I_3 \leq    \int_{B_r(x_0)}   |\nabla u|^{p(x)-1}  .|\nabla v| dx
~\leq~  {1\over 4}  \int_{B_r(x_0)}  |\nabla u|^{p(x)} dx + c\int_{B_r(x_0)}  |\nabla v|^{p(x)}dx.
\end{align}
\n Using (2.5)-(2.8), we get
\begin{eqnarray}\label{2.9}
\int_{B_r(x_0)} |\nabla u|^{p(x)} dx&\leq& 4c \int_{B_r(x_0)}
|\nabla v|^{p(x)} dx-2\int_{B_{R}(x_0)} g.(u-v)dx\nonumber\\
&&~~+2\int_{B_{R}(x_0)}  F.\nabla ( u-v)dx.
\end{eqnarray}
\n Combining (2.2), and (2.9), we get for a positive constant
\[\displaystyle{C_1=C_1(n,s,p_-,p_+,\text{diam}(\Omega),||\nabla p||_{L^s(\Omega)},||\nabla u||_{L^{p(x)}(\Omega)})}\]
\begin{eqnarray}\label{2.10}
\int_{B_r(x_0)} |\nabla u|^{p(x)} dx&\leq& C_1\Big(\Big( \Big({ r\over R}\Big)^n +R^{\beta\over 2}\Big)\int_{B_R(x_0)} |\nabla u|^{p(x)}dx +  R^{n+\beta}\Big)\nonumber\\
&&~-2\int_{B_{R}(x_0)} g.(u-v)dx+2\int_{B_{R}(x_0)}  F.\nabla ( u-v)dx.
\end{eqnarray}

\n Applying H\"{o}lder's inequality and using (1.6) and the fact that
$|u|, |v|\leq M$ in $B_{R}(x_0)$, we obtain
\begin{eqnarray}\label{2.11}
&&\left|-2\int_{B_{R}(x_0)} g.(u-v)dx\right| \leq 4M \int_{B_{R}(x_0)} |g|dx\nonumber\\
&&\qquad\leq 4M ||g||_{L^{t_1}(B_{R}(x_0))}.|B_{R}(x_0)|^{1-{1\over{t_1}}}\nonumber\\
&&\qquad=4M ||g||_{L^{t_1}(B_{R}(x_0))}.\omega_n^{1-{1\over{t_1}}}.R^{n-{n\over{t_1}}}\nonumber\\
&&\qquad =C_2R^{n-p_m+\alpha_1 p_m},\qquad \alpha_1=1-{{n}\over{t_1p_m}},\quad C_2=4M \omega_n^{1-{1\over{t_1}}}.||g||_{L^{t_1}(\Omega)}
 \end{eqnarray}

\n Due to (1.5) and the continuity of $p(x)$, we have $\displaystyle{t_1>{{n}\over{p_m}} }$.
Hence $\alpha_1\in(0,1)$.

\n We observe from (1.4) and (1.7) that we have
\[t_2\geq{{t_2p(x)}\over{n}}>{{p(x)}\over{p(x)-1}}\geq{{p_m}\over{p_m-1}}\]

\n Applying Young's and H\"{o}lder's inequalities and taking into account (1.8)
and using the convexity of $|\xi|^{p_m}$, we obtain for $\epsilon\in(0,1)$ and a constant $C_3$ depending
only on $p_-$ and $p_+$
\begin{eqnarray}\label{2.12}
&&\left|2\int_{B_{R}(x_0)} F.\nabla ( u-v)dx\right| \leq C_3 \int_{B_{R}(x_0)}|F|^{{p_m}\over{p_m-1}}dx
+\epsilon \int_{B_{R}(x_0)}|\nabla(u-v)|^{p_m}  dx\nonumber\\
&&\quad\leq C_3||F||_{L^{t_2}(\Omega)}^{{p_m}\over{p_m-1}}.|B_R(x_0)|^{{t_2(p_m-1)-p_m}\over{t_2(p_m-1)}}+\epsilon\int_{B_{R}(x_0)}|\nabla(u-v)|^{p_m}dx\nonumber\\
&&\quad=C_3\omega_n^{{t_2(p_m-1)-p_m}\over{t_2(p_m-1)}}||F||_{L^{t_2}(\Omega)}^{{p_m}\over{p_m-1}} R^{n-{{np_m}\over{t_2(p_m-1)}}}+\epsilon \int_{B_{R}(x_0)}|\nabla(u-v)|^{p_m}dx\nonumber\\
&&\quad\leq C_3' R^{n-p_m+\alpha_2 p_m} +
 \epsilon 2^{p_m-1} \int_{B_{R}(x_0)}|\nabla u|^{p_m}dx+\epsilon 2^{p_m-1} \int_{B_{R}(x_0)}|\nabla v|^{p_m}dx
 \end{eqnarray}

\n where $\displaystyle{\alpha_2=1-{{n}\over{t_2(p_m-1)}} }$ and $\displaystyle{C_3'=C_3.\omega_n^{{t_2(p_m-1)-p_m}\over{t_2(p_m-1)}}.||F||_{L^{t_2}(\Omega)}^{{p_m}\over{p_m-1}} }$.

\n Due to (1.7), and the continuity of $p(x)$, we have $\displaystyle{t_2>{{n}\over{p_m-1}} }$,
and therefore $\alpha_2\in(0,1)$

\vs 0.2cm\n Observe that we have for $w\in W^{1,p(x)}(B_{R}(x_0))$
\begin{eqnarray}\label{2.13}
\int_{B_{R}(x_0)}|\nabla w|^{p_m}dx&=& \int_{B_{R}(x_0)\cap[|\nabla w|\leq 1]}|\nabla w|^{p_m}dx
+\int_{B_{R}(x_0)\cap[|\nabla w|>1]}|\nabla w|^{p_m}dx\nonumber\\
   &\leq & \int_{B_{R}(x_0)}|\nabla w|^{p(x)}  dx+\omega_n R^n.
 \end{eqnarray}

\n Using (2.3) and (2.13), we get from (2.12)
\begin{eqnarray}\label{2.14}
&&\left|2\int_{B_{R}(x_0)} F.\nabla ( u-v)dx\right| ~\leq C_3' R^{n-p_m+\alpha_2 p_m}+\epsilon 2^{p_m-1}(1+{p_+\over p_-}) \int_{B_{R}(x_0)}|\nabla u|^{p(x)}  dx+\epsilon2^{p_m}\omega_nR^n\nonumber\\
&&\qquad=(C_3'+\epsilon2^{p_m}\omega_nR^{(1-\alpha_2)p_m}) R^{n-p_m+\alpha_2 p_m}+\epsilon 2^{p_m-1}(1+{p_+\over p_-}) \int_{B_{R}(x_0)}|\nabla u|^{p(x)}dx\nonumber\\
&&\qquad\leq C_4\Big(R^{n-p_m+\alpha_2 p_m}+\epsilon \int_{B_{R}(x_0)}|\nabla u|^{p(x)}dx\Big)
\end{eqnarray}

\n where $\displaystyle{C_4=\max\left(C_3'+2^{p_m}\omega_n(\text{diam}(\Omega)/2)^{(1-\alpha_2)p_m},2^{p_m-1}(1+{p_+\over p_-})\right)}$.

\n Combing (2.10), (2.11) and (2.14), we get
\begin{eqnarray}\label{2.15}
&&\int_{B_r(x_0)} |\nabla u|^{p(x)} dx~\leq~ C_1\Big(\Big( \Big({ r\over R}\Big)^n +R^{\beta\over 2}\Big)\int_{B_R(x_0)} |\nabla u|^{p(x)}dx\Big)\nonumber\\
&&\qquad+C_1R^{n+\beta}+C_2R^{n-p_m+\alpha_1 p_m}+C_4R^{n-p_m+\alpha_2 p_m}+C_4\epsilon \int_{B_{R}(x_0)}|\nabla u|^{p(x)}dx
\end{eqnarray}

\n Setting $\phi(r)=\displaystyle{\int_{B_r(x_0)} |\nabla u|^{p(x)}dx}$,
we obtain from (2.15)
\begin{eqnarray}\label{2.16}
\phi(r)\leq \max(C_1,C_2,C_4)\Big(\Big( \Big({ r\over R}\Big)^n +R^{\beta\over 2}+\epsilon\Big)\phi(R)
+R^{n+\beta}+R^{n-p_m+\alpha_1 p_m}+R^{n-p_m+\alpha_2 p_m}\Big).\nonumber\\
&&
\end{eqnarray}

\n Now observe that we have for $\alpha=\min(\alpha_1,\alpha_2)$
\begin{eqnarray}\label{2.17}
&&R^{n+\beta}+R^{n-p_m+\alpha_1 p_m}+R^{n-p_m+\alpha_2 p_m}=
(R^{\beta+(1-\alpha)p_m}+R^{(\alpha_1-\alpha)p_m}+R^{(\alpha_2-\alpha)p_m})R^{n-p_m+\alpha p_m}\nonumber\\
&&~~\leq \Big((\text{diam}(\Omega)/2)^{\beta+(1-\alpha)p_m}+(\text{diam}(\Omega)/2)^{(\alpha_1-\alpha)p_m}
+(\text{diam}(\Omega)/2)^{(\alpha_2-\alpha)p_m}\Big)R^{n-p_m+\alpha p_m}\nonumber\\
&&~~=C_4'R^{n-p_m+\alpha p_m}
\end{eqnarray}

\n Using (2.16) and (2.17), we get for $C_5=\max(C_1,C_2,C_4).\max(1,C_4')$
\begin{equation}\label{2.18}
\phi(r) \leq C_5\Big(\Big({r\over R}\Big)^n+R^{\beta\over 2}+\epsilon\Big) \phi(R)+ R^{n-p_m+\alpha p_m}\Big)
\quad\forall r\in(0,R).
\end{equation}
\n Given that $\phi$ is a nonnegative and a nondecreasing function on
$(0,R)$, we infer from (2.18) and Lemma 5.12 p. 248 of \cite{[MZ]}, that we have
for a constant $C_6>0$
\begin{equation}\label{2.19}
\phi(r) \leq C_6\Big({r\over R}\Big)^{n-p_m+\alpha
p_m}\big(\phi(R) + R^{n-p_m+\alpha p_m} \big) \qquad \forall r\in(0,R).
\end{equation}

\n Using (2.19) and (2.13), for $R=r$ and $w=u$, we get for
$\displaystyle{C_7={{C_6(\phi(R) + R^{n-p_m+\alpha p_m})}\over{R^{n-p_m+\alpha
p_m}}} }$
\begin{eqnarray*}
\int_{B_r(x_0)} |\nabla u|^{p_m} dx &\leq & \int_{B_r(x_0)} |\nabla u|^{p(x)} dx+\omega_nr^n
~\leq~ C_7r^{n-p_m+\alpha p_m}+\omega_nr^n\nonumber\\
&=& (C_7+\omega_nr^{(1-\alpha)p_m})r^{n-p_m+\alpha p_m}~\leq~ Cr^{n-p_m+\alpha p_m}
\end{eqnarray*}
\n where $C=C_7+\omega_n(\text{diam}(\Omega)/2)^{(1-\alpha)p_m}$. This completes the proof of the lemma.
\qed

\vs 0,3cm\n\emph{Proof of Theorem 1.1.}
Using (2.1) and H\"{o}lder's inequality, we obtain for all $R\in(0,R_1)$
\begin{eqnarray*}
\forall r\in(0,R)~~~\int_{B_r}|\nabla u|dx &\leq&|B_r|^{1-{1\over {p_m}}}\Bigg(\int_{B_r}|\nabla u|^{p_m}dx\Bigg)^{1\over {p_m}}\nonumber\\
&\leq& \omega_n^{1-{1\over {p_m}}}.r^{n-{n\over {p_m}}}.(Cr^{n-p_m+\alpha p_m})^{1\over{p_m}}\nonumber\\
&=& C^{1\over{p_m}}\omega_n^{1-{1\over {p_m}}}.r^{n-{n\over {p_m}}}.r^{{n-p_m+\alpha p_m}\over{p_m}}\nonumber\\
&\leq& C'r^{n-1+\alpha},\qquad C'=C^{1\over{p_m}}\omega_n^{1-{1\over {p_m}}}.
\end{eqnarray*}

\n We conclude (\cite{[MZ]} Theorem(Morrey) 1.53 p. 30) that $u\in
C^{0,\alpha}_{loc}(\Omega)$, which completes the proof of the theorem.\qed

\vs 0.3cm

\begin{remark}
If $g, F\in L^{\infty}(\Omega)$, then (1.5)-(1.8)
are satisfied for any $t_1>{{n}\over{p(x)}}$ and $t_2>{{n}\over{p(x)-1}}$. Therefore,
we obtain $u\in C^{0,\alpha}_{loc}(B_R(x_0))$ for any $B_R(x_0)\subset\subset\Omega$ with
$\displaystyle{\alpha=1-\max\left({{n}\over{t_1p_m}},{{n}\over{t_2(p_m-1)}}\right) }$,
$\displaystyle{p_m=\min_{x\in \overline{B}_R(x_0)}p(x)}$ and $R$ small enough. Given that
$t_1$ and $t_2$ can be chosen arbitrarily large, we obtain
$u\in C^{0,\alpha}_{loc}(\Omega)$ for any $0<\alpha<1$.
\end{remark}

\end{document}